\documentclass[11pt]{article}
\usepackage{amsmath,amssymb}
\newcommand{\mind}{\mathrm{Ind\,}}

\newcommand{\grad}{{\rm \nabla}}

\newcommand{\real}{\mathbb{R}}

\newcommand {\cqd}{\begin{flushright}\vskip-25pt$\Box$\end{flushright}}
\newtheorem{myth}{Theorem}[section]
\newtheorem{mylem}{Lemma}[section]
\newtheorem{myprop}{Proposition}[section]
\newtheorem{mydef}{Definition}[section]
\newtheorem{myrem}{Remark}[section]
\newtheorem{mycoro}{Corollary}[section]
\newtheorem{myquest}{Question}[section]
\def\aco{\lbrace }
\def\acf{\rbrace }
\begin{document}
\title{ Stable hypersurfaces with constant scalar curvature in Euclidean spaces}\author { Hil\'{a}rio Alencar\thanks{The authors were partially supported by CNPq and FAPERJ, Brazil.\newline \indent MSC2010 Classification: 53C42},  Walcy
Santos$^*$ {and} {Detang Zhou$^*$ }}
\date {\today}
\maketitle \vspace*{7pt}
\begin {quote}
{\scriptsize{\bf Abstract. } We obtain some nonexistence results for
complete noncompact  stable hyppersurfaces with nonnegative constant
scalar curvature in Euclidean spaces. As a special case we prove
that there is no complete noncompact strongly stable hypersurface
$M$ in $\mathbb{R}^{4}$ with zero scalar curvature $S_2$, nonzero
Gauss-Kronecker curvature and finite total curvature (i.e.
$\int_M|A|^3<+\infty$). }

\end{quote}
\vspace*{7pt}
\section {Introduction}\indent
In this paper we study the complete noncompact stable hypersurfaces
with constant scalar curvature in  Euclidean spaces. It has been
proved by Cheng and Yau \cite{CY} that any complete noncompact
hypersurfaces in the Euclidean space with constant scalar curvature
and nonnegative sectional curvature must be a generalized cylinder.
Note that the assumption of nonnegative sectional curvature is  a
strong condition for hypersurfaces in the Euclidean space with zero
scalar curvature since the hypersurface has to be flat in this case.
Let $M^n$ be a complete orientable Riemannian manifold and let
$x:M^n\rightarrow\mathbb{R}^{n+1}$ be an isometric immersion into
the Euclidean space $\mathbb{R}^{n+1}$ with constant scalar
curvature.  We can choose a a global unit normal vector field $N$ and the
Riemannian connections  $\nabla$ and $\widetilde{\nabla}$ of $M$ and
$\mathbb{R}^{n+1}$, respectively, are related by
$$\widetilde{\nabla}_XY=\nabla_XY+\langle A(X),Y\rangle N,$$
where $A$ is the second fundamental form of the immersion, defined
by
$$A(X)=-\widetilde{\nabla}_XN.$$

Let ${\lambda _1,...,\lambda _n }$ be the eigenvalues of $A$. The
{\it r-mean curvature} of the immersion in a point $p$ is defined by
$$H_r=\frac{1}{\binom nr}\sum_{i_1<...<i_r}\lambda _{i_1}...\lambda
_{i_r}=\frac{1}{\binom nr}S_r,$$ where $S_r$ is the $r$-symmetric
function of the ${\lambda _1,...,\lambda _n }$, $H_0=1$ and $H_r=0$,
for all $r\geq n+1$. For $r=1$, $H_1=H$ is the mean curvature of the
immersion, in the case $r=2$, $H_2$ is the normalized scalar
curvature and for $r=n$, $H_n$ is the Gauss-Kronecker curvature.

It is well-known  that hypersurfaces with constant scalar
curvature are critical points for a geometric variational problem,
namely, that associated to the functional
\begin{equation}\label{rarea}
\mathcal{A}_1(M)=\int_{M}S_1\,dM,
\end{equation}
under compactly supported variations that preserves volume.  Let $M$
be a hypersurface in the Euclidean space  with constant scalar
curvature. Following \cite{AdCE},  when the scalar curvature is
zero, we say that a regular domain $D\subset M$ is {\it stable} if
the critical point is such that
$(\frac{d^2\mathcal{A}_1}{dt^2})_{t=0}\ge 0$, for all variations
with compact support in $D$ and when the scalar curvature is
nonzero, we say that a regular domain $D\subset M$ is {\it strongly
stable } if the critical point is such that
$(\frac{d^2\mathcal{A}_1}{dt^2})_{t=0}\ge 0$, for all variations
with compact support in $D$. It is natural to study the global
properties of hypersurfaces in  the Euclidean space with constant
scalar curvature. For example we have the following open question
(see 4.3 in \cite{AdCE}).
\begin{myquest}\label{quest1.1}Is there  any   stable complete hypersurfaces $M$ in $\real^{4}$ with zero scalar curvature and nonzero the Gauss-Kronecker
curvature?
\end{myquest}
We have a partial answer to the question \ref{quest1.1}.\\[5pt]

\noindent {\bf Theorem A.} (see Theorem \ref{thm3.1}) {\em There is
no complete noncompact stable hypersurface $M$ in $\mathbb{R}^{n+1}$
with zero scalar curvature $S_2$ and $3$-mean curvature $S_3\ne 0$
satisfying
\begin{equation}
\lim_{R\to +\infty}\frac{\int_{B_R}S_1^3}{R^2}=0,
\end{equation}
where $B_R$ is the geodesic ball in $M$.}
\\[5pt]
When $S_2=0$, $S_1^2=|A|^2$ we have\\[5pt]
\noindent {\bf Corollary B.}
{\em There is no complete noncompact stable hypersurface $M$ in
$\mathbb{R}^{4}$ with zero scalar curvature $S_2$, nonzero
Gauss-Kronecker curvature and finite total curvature (i.e.
$\int_M|A|^3<+\infty$)}.
\\[5pt]
We remark that Shen and Zhu (see \cite{SZ}) proved that a complete stable minimal $n$-dimensional hypersurface in $\mathbb{R}^{n+1}$ with finite total curvature must be a hyperplane. The above Corollary can be seen as a similar result in dimension $3$ for hypersurfaces with zero scalar curvature.

 We also prove the following result for hypersurfaces with positive constant scalar curvature in Euclidean space.\\[5pt]
\noindent {\bf Theorem C.} (see Theorem \ref{thm3.2})  {\em There is
no complete immersed strongly stable hypersurface $M^n\rightarrow
\mathbb{R}^{n+1}$, $n\geq 3$, with positive constant scalar
curvature and polynomial growth of $1$-volume, that is
$$\lim_{R\rightarrow\infty}\frac{\int _{B_R}S_1dM}
{R^n}<\infty,$$ where $B_R$ is a geodesic ball of radius $R$ of
$M^n$.}\\[5pt]

As a consequence of the properties of a graph with constant scalar curvature, we have the following corollary:\\[5pt]
\noindent {\bf Corollary D.} (see Corollary \ref{coroD})
{\em Any entire
graph on $\mathbb{R}^n$ with nonnegative constant scalar curvature
must have zero scalar curvature.}
\\[5pt]
This can be compared with a result
of Chern \cite{Ch} which says any entire graph on $\mathbb{R}^n$
with constant mean curvature must be minimal. It has been be known
by a result of X. Cheng in \cite{Che} (see also \cite{ENR}) that  any
complete noncompact stable hypersurface in $\mathbb{R}^{n+1}$ with
constant mean curvature must be minimal if $n<5$. It is natural to
ask that any complete noncompact stable hypersurface in
$\mathbb{R}^{n+1}$ with nonnegative constant scalar curvature must
have zero scalar curvature.

It should be remarked that Chern\cite{Ch}  proved that there is no  entire graph on $\mathbb{R}^n$  with  Ricci  curvature less than a negative constant. We don't know whether there exists an entire graph on $\mathbb{R}^n$  with constant negative scalar curvature.

The rest of this paper is organized as follows: we include some
results and definitions which will be used in the proof of our
theorems in Section \ref{secind}. The proof of main results are
given in Section 3 and Section 4 is an appendix in which we prove
some stability properties for graphs with constant scalar curvature
in the Eucildean space.

\section{Some stability and index properties for hypersurfaces with $S_2=const$.}\label{secind}

We introduce the $r$'th Newton transformation, $P_r:T_pM\rightarrow
T_pM$, which are defined inductively by
$$\begin{array}{ll}
P_0= & I,\\
P_r= &  S_rI-A\circ P_{r-1},\;r\geq 1.
\end{array}$$
The following formulas are useful in the proof (see, \cite{Re}, Lemma
2.1).
\begin{eqnarray}
    \;\mbox{{\rm
trace}}(P_r)&=&(n-r)S_r,\label{eq0.1}  \\
  \;\mbox{{\rm
trace}}(A\circ P_r)&=&(r+1)S_{r+1}, \label{eq0.2} \\
\;\mbox{{\rm trace}}(A^2\circ P_r)&=&S_1S_{r+1}-(r+2)S_{r+2}.
\label{eq0.3}
\end{eqnarray}

 From \cite{AdCC} we have
 the second variation formula for hypersurfaces in a space form of constant curvature $c$,
 $\mathbb{Q}_c^{n+1}$,
with constant $2$-mean curvature:
\begin{equation}\label{secvar}
   \frac{d^2\mathcal{A}_1}{dt^2}|_{t=0}= \int_D\langle P_1(\grad f),\grad f\rangle dM- \int_D(S_1S_2-3S_3+c(n-1)S_1)f^2dM, \quad \forall f\in C_c^\infty(D).
\end{equation}

\begin{mydef} When $S_2=0$ and $c=0$, $M$ is stable if and only if
\begin{equation}\label{eqn15}
    \int_M\langle P_1(\nabla f), \nabla f\rangle dM\ge -3\int_MS_3f^2dM,
\end{equation}
for any $f\in C_c^{\infty}(M)$. One can see  that if $P_1\equiv 0$,
then $S_3=0$ and $M$ is
 stable.
When $S_2=const.\neq 0$, $M$ is  stable if and only if
$$
    \int_D\langle P_1(\grad f),\grad f\rangle dM\geq \int_D(S_1S_2-3S_3+c(n-1)S_1)f^2dM,
$$
for all $f\in C_c^{\infty}(M)$ and $\int_M f dM=0$.  We say that $M$
is strongly stable if and only if the above inequality holds for all
$f\in C_c^{\infty}(M)$.
\end{mydef}
Similar to minimal hypersurface we can also define the index $I$ for
hypersurfaces with constant scalar curvature. Given a relatively
compact domain $\Omega \subset M$, we denote by $\mind^1(\Omega)$
the number of linearly independent normal deformations with support
on $\Omega$ that decrease $\mathcal{A}_1$.  The \emph{index} of the
immersion are defined as
\begin{equation}\label{E:swi}
\mind^1(M) := \sup \aco \mind^1(\Omega) ~\big|~ \Omega
        \subset M, ~~\Omega \mbox{~~relatively compact}\acf . \\
\end{equation}
$M$ is \emph{strongly stable}  if $\mind^1(M)=0$. The following
result has been known in  \cite{El}.
\begin{mylem} Let $M^n\rightarrow\mathbb{Q}_c^{n+1}$ be a noncompact hypersurface with $S_2=const.>0$.
 If $M$ has finite index then there exist a compact set $K\subset M$ such that $M\setminus K$
 is strongly stable.\label{lemma2}
\end{mylem}
For hypersurfaces with constant mean
curvature,  do Carmo and Zhou \cite{dCZ} proved that
\begin{myth} { Let $x: M^n \rightarrow \overline{M}^{n+1}$
 be an isometric immersion with constant mean curvature $H$.
 Assume $M$ has subexponential volume growth and finite index. Then there
exist a constant $R_0$ such that   }
$$H\leq - \overline{Ric}_{M\setminus B_{R_0}}(N),$$
{ where $N$ is a smooth normal vector field along $M$ and
$\overline{Ric}(N)$ is the Ricci curvature of $\overline{M}$ in the
normal vector $N$.}
\end{myth}
The technique in \cite{dCZ} was generalized by Elbert \cite{El} to
prove  the following result:
\begin{myth}\label{thm2.2} { Let $x: M^n \rightarrow \mathcal{Q}(c)^{n+1}$
 be an isometric immersion with
$S_2 = constant > 0$. Assume that $\mind^1M < \infty$ and that the
1-volume of M is infinite and has polynomial growth. Then c is
negative and }
$$S_2\leq -c.$$
\end{myth}

In particular, it implies that when $c=0$  the hypersurfaces in the
above theorem must have nonpositive scalar curvature.

\section{Proof of  the theorems}

When $S_2=0$ we know that $|S_1|^2=|A|^2$. Thus, if $S_3\ne 0$,
we have that $|A|^2>0$. Hence $S_1\neq 0$ and we can choose an
orientation such that  $P_1$ is semi-positive definite. Since
\begin{equation*}
    \begin{split}
|\sqrt{P_1}A|^2&=\mbox{{\rm
trace}}(A^2\circ P_1)\\
&=-3S_3,
     \end{split}
\end{equation*} then, when $c=0$,  $M$ is stable if
\begin{equation}\label{est}
    \int_M\langle P_1(\nabla f), \nabla f\rangle
    dM\ge\int_M|\sqrt{P_1}A|^2f^2dM,
\end{equation}
for any $f\in C_c^{\infty}(M)$.

 By Lemma 4.1 in \cite{AdCC}, when $S_2=0$, we know that $|\nabla A|^2-|\nabla S_1|^2\ge 0$. In the following lemma, we  characterize the equality case in some special case.

\begin{mylem}\label{lem1} Let $M^n(n\ge 3)$ be a non-flat  connected immersed 1-minimal hypersurface in
$\real^{n+1}$. If $|\nabla A|^2=|\nabla S_1|^2 $ holds  on all
nonvanishing point of $|A|$ in $M$, then  each component of $M$ with
$|A|\ne 0$ must be a  cylinder over a curve.
\end{mylem}
 {\bf Proof.}\ Choose a
frame at $p$ so that the second fundamental form is diagonalized. From the computations in \cite{SSY}, we have
$|A|^2=\sum_{i}h_{ii}^2$, and
\begin{eqnarray}\label{eqnlong}
   \sum_{i,j,k}h_{ijk}^2 -|\nabla|A||^2&=& [(\sum_{i,j}h_{ij}^2)(\sum_{s,t,k}h_{stk}^2)-\sum_k(\sum_{i,j}h_{ij}h_{ijk})^2](\sum_{i,j}h_{ij}^2)^{-1} \nonumber\\
   &=& \frac 12\sum_{i,j,k,s,t}(h_{ij}h_{stk}-h_{st}h_{ijk})^2 |A|^{-2}\nonumber\\
   &=&\frac 12\left[\sum_{i,k,s,t}(h_{ii}h_{stk}-h_{st}h_{iik})^2
   +\sum_{s}h_{ss}^2(\sum_{k}\sum_{i\ne j}h_{ijk}^2)\right]|A|^{-2}\nonumber\\
   &=&\frac 12\left[\sum_{i,k,s}(h_{ii}h_{ssk}-h_{ss}h_{iik})^2+\sum_{i}h_{ii}^2(\sum_{k}\sum_{s\ne t}h_{stk}^2)\right]|A|^{-2}\nonumber\\
   &&\qquad\qquad+\frac 12
   (\sum_{k}\sum_{i\ne j}h_{ijk}^2)   \nonumber\\
   &=&\frac
   12\left[\sum_{i,k,s}(h_{ii}h_{ssk}-h_{ss}h_{iik})^2\right]|A|^{-2}+(\sum_{k}\sum_{i\ne
   j}h_{ijk}^2)\nonumber\\
   &=&\frac
   12\left[\sum_{i,k,s}(h_{ii}h_{ssk}-h_{ss}h_{iik})^2\right]|A|^{-2}+\nonumber\\
   &&2\sum_{i\ne
   j}h_{iij}^2+\sum_{i\ne
   j,j\ne k, i\ne k}h_{ijk}^2\ge 0.
\end{eqnarray}
It is clear that  the right hand side is nonnegative and is zero if
and only if all terms on the right hand side vanish.
\begin{equation}\label{eqn19}
    \sum_{i,j,k}h_{ijk}^2 -|\nabla|A||^2\ge 0.
\end{equation}

  Suppose $x: M\to \real^{n+1}$ is the 1-minimal
immersion. Since $M$ is not a hyperplane, then $|A|$ is a
nonnegative continuous function which does not vanish identically.
Let $p$ be such a point such that $|A|(p)>0$. Then $|A|>0$ in a
connected open set $U$ containing $p$.
 The equality in (\ref{eqn19}) implies
\begin{equation*}\begin{split}
    h_{jji}&=0, \textrm{ for all }j\ne i,\\
    h_{ijk}&=0, \textrm{ for all }j\ne i,j\ne k, k\ne i\\
    h_{ii}h_{ssk}&=h_{ss}h_{iik}, \textrm{ for all }i,s,k.
    \end{split}
\end{equation*}
 So we have $ h_{jij} =0,$   for all  $j\ne i,$ and from the last
 equation we claim that at most one $i$ such that $h_{iii}\ne 0$.
 Otherwise, without the loss of generality we assume  $h_{111}\ne 0$, and $h_{222}\ne 0$, we have $h_{11}h_{22k}=h_{22}h_{11k}$ for all $k$. This implies $h_{11}=h_{22}=0$ by choosing $k=1,2$. Using again the third formula we have $h_{jj}h_{111}=h_{11}h_{jj1}$ for $j=3, \cdots, n$. Hence $h_{jj}=0 $ for all $j=3, \cdots, n$, which contradicts to $|A|\ne 0$.

 We now assume $h_{111}\ne 0$ by continuity we can also assume $h_{11}\ne 0$.  From the last
 equation of above equation, we have $h_{11}h_{ss1}=h_{ss}h_{111}$ for $s\ne
 1$. Hence $h_{ss}=0$ for all $s\ne 1$. This implies that $M$ is a   cylinder
over a curve.
\\ \cqd
 We are now ready to prove

\begin{myth}\label{thm3.1}
There is no complete noncompact stable hypersurfaces in $R^{n+1}$
with $S_2=0$ and $S_3\ne 0$ satisfying
\begin{equation*}
\lim_{R\to +\infty}\frac{\int_{B_R}S_1^3}{R^2}=0.
\end{equation*}
 \end{myth}
{\bf Proof.} Assume for the sake of contradiction that there were
such a hypersurface $M$. From Lemma 3.7 in \cite{AdCC}, we have
\begin{equation}\label{eqn16}
    L_1S_1=|\nabla A|^2-|\nabla S_1|^2+3S_1S_3.
\end{equation}
Since for any $\phi\in C_c^{\infty}(M)$,
\begin{equation*}
    \begin{split}
       \int_M\langle P_1(\nabla (\phi S_1)),\nabla (\phi S_1)\rangle dM  & =\int_M\langle P_1((\nabla \phi )S_1+\phi \nabla S_1),(\nabla\phi )S_1)+\phi \nabla S_1 \rangle dM\\
         & =\int_M\langle P_1(\nabla \phi ),\nabla \phi \rangle
         S_1^2 dM+2 \int_M\langle P_1(\nabla \phi ),\nabla S_1\rangle
         \phi S_1 dM\\
                  &\qquad +\int_M\phi^2 \langle P_1(\nabla S_1),\nabla
                  S_1\rangle dM,
              \end{split}
\end{equation*}
then using (\ref{eqn16}) we have
\begin{equation*}
    \begin{split}
       \int_M\phi^2 & S_1(|\nabla A|^2-|\nabla S_1|^2) dM  = \int_M (L_1S_1-3S_1S_3)\phi^2 S_1 dM\\
         & =-\int_M\langle P_1(\nabla S_1),\nabla (\phi^2S_1 )\rangle dM
         - \int_M3S_3\phi^2  S_1^2 dM\\
                  &=-\int_M\phi^2 \langle P_1(\nabla S_1),\nabla
                  S_1\rangle dM-2 \int_M\langle P_1(\nabla \phi ),\nabla S_1\rangle
         \phi S_1 dM- \int_M3S_3\phi^2  S_1^2 dM\\
         &=-\int_M\langle P_1(\nabla (\phi S_1)),\nabla (\phi S_1)\rangle dM+\int_M\langle P_1(\nabla \phi ),\nabla (\phi )\rangle
         S_1^2 dM- \int_M3S_3\phi^2  S_1^2 dM\\
         &\le \int_M\langle P_1(\nabla \phi ),\nabla \phi \rangle
         S_1^2 dM\\
         & \le \int_M|\nabla \phi |^2
         S_1^3 dM,
              \end{split}
\end{equation*}
 for any $\phi\in C_c^{\infty}(M)$.  Here we have
 used the stability inequality (\ref{eqn15}) in the fifth line and use the following consequence of (\ref{eq0.1}) in the last inequality:

 \begin{equation}\label{estima}
             S_1|\nabla \phi|^2\ge\langle P_1(\nabla \phi), \nabla
             \phi\rangle.
 \end{equation}
 We can choose $\phi$ as
\begin{equation*}
    \phi(x)=\left\{
                 \begin{array}{ll}
                   \frac{2R-r(x)}{R}, & \hbox{on } B_{2R}\setminus B_R; \\
                   1, & \hbox{on } B_R; \\
                   0, & \hbox{on }M\setminus B_{2R}.
                 \end{array}
               \right.
\end{equation*}
Thus from the choice of $\phi$ we have $S_1(|\nabla A|^2-|\nabla
S_1|^2)\equiv 0$. Therefore the elipticity of $L_1$ implies
$L_1S_1=3S_1S_3$. From Lemma \ref{lem1}, $M$ must be a  cylinder
over a curve which contradicts $S_3\ne 0$. The proof is complete. \\
\cqd

The following Lemma is of some independent interest and we include
here since its second part is useful in the proof of Theorem
\ref{thm3.2}.

\begin{mylem}\label{lem3.2}
Let $M$ be a complete immersed hypersurface in $\mathcal{Q}_c^{n+1}$
with nonnegative constant scalar curvature $S_2>-\frac{n(n-1)}2c$
and $S_1\ne 0$.  \newline 1)If M is strongly stable outside a
compact subset, then either $M$ has finite 1-volume,
or\begin{equation*}
    \lim_{R\to +\infty}\frac{1}{R^2}\int_{B_{R}}S_1=+\infty.
\end{equation*}
 2)If M is  strongly stable,
then \begin{equation*}
    \lim_{R\to +\infty}\frac{1}{R^2}\int_{B_{R}}S_1=+\infty.
\end{equation*}In particular $M$ has infinite 1-volume.

\end{mylem}
{\bf Proof.} We can assume that there exists a geodesic ball $B_{R_0}\subset M$ such
that $M\setminus B_{R_0}$ is strongly stable. That is,
\begin{equation}\label{est9}
\int_M (S_1S_2-3S_3+c(n-1)S_1)f^2dM \leq
 \int_M\langle P_1(\grad f),\grad f\rangle dM,
\end{equation}
for all $f\in \mathcal{C}_c^\infty(M\setminus B_{R_0})$.

Now, since $S_2\ge 0$, we have (see \cite{AdCR}, p. 392)
$$H_1H_2\geq H_3,$$
and $$H_1\geq H_2^{1/2}.$$
 By using that $S_1=nH_1$,
$S_2=\dfrac{n(n-1)}{2}H_2$ and $S_3=\dfrac{n(n-1(n-2)}{6}H_3$, it
follows that
$$ \frac{(n-2)}{n}S_1S_2\geq 3S_3,$$
that is,
\begin{equation}\label{est11}
    -3S_3\geq-\frac{(n-2)}{n}S_1S_2.
\end{equation}
We also have that
$$\frac{S_1}{n}\geq\left( \frac {2S_2}{n(n-1)}\right)^{1/2},$$
which implies
\begin{equation}\label{est12}
    S_1\geq \left(\frac{2n}{n-1}  \right)^{1/2}S_2^{1/2}.
\end{equation}
By using inequality (\ref{est11}) in (\ref{est9}), it follows that
$$ \int_M \left(S_1S_2-\frac{n-2}{n}S_1S_2+c(n-1)S_1\right)f^2dM\leq
\int_M\langle P(\nabla f),\nabla f\rangle dM,$$
that is,
$$ \int_M \left(S_2+\frac{n(n-1)c}{2}\right)S_1f^2dM\leq\frac{n}{2}
\int_M\langle P(\nabla f),\nabla f\rangle dM.$$ By using
(\ref{estima}), we obtain that $$ \int_MS_1|\nabla f|^2dM\ge
\int_M\langle P(\nabla f),\nabla f\rangle dM$$
 Therefore,
there exists a constant $C>0$ such that
\begin{equation}\label{eq}
     \int_MS_1|\nabla f|^2dM\ge C\int_MS_1f^2dM.
 \end{equation}
1) When $M$ is strongly stable outside $B_{R_0}$. We can choose $f$
as
\begin{equation*}
    f(x)=\left\{
                 \begin{array}{ll}
                   r(x)-R_0, & \hbox{on } B_{R_0+1}\setminus B_{R_0}; \\
                   1, & \hbox{on } B_{R+R_0+1}\setminus B_{R_0+1}; \\
                   \frac{2R+R_0+1-r(x)}{R}, & \hbox{on } B_{2R+R_0+1}\setminus B_{R+R_0+1}; \\
                   0, & \hbox{on }M\setminus B_{2R+R_0+1},
                 \end{array}
               \right.
\end{equation*}
where $r(x)$ is the distance function to a fixed point. Then

\begin{equation*}
\frac{1}{R^2}\int_{B_{2R+R_0+1}\setminus
B_{R+R_0+1}}S_1dM+\int_{B_{R_0+1}\setminus B_{R_0}}S_1dM\ge
C\int_{B_{R+R_0+1}\setminus B_{R_0+1}}S_1dM.
\end{equation*}
If the 1-volume is infinite, we can choose $R$ large such that
$$C\int_{B_{R+R_0+1}\setminus B_{R_0+1}}S_1dM>\int_{B_{R_0+1}\setminus
B_{R_0}}S_1dM,$$ hence
\begin{equation*}
    \lim_{R\to +\infty}\frac{1}{R^2}\int_{B_{2R+R_0+1}\setminus B_{R+R_0+1}}S_1dM=+\infty.
\end{equation*}

2) When $M$ is strongly stable we can choose a simpler test function
$f$ as
\begin{equation*}
    f(x)=\left\{
                 \begin{array}{ll}
                   1, & \hbox{on } B_{R}; \\
                   \frac{2R-r(x)}{R}, & \hbox{on } B_{2R}\setminus B_{R}; \\
                   0, & \hbox{on }M\setminus B_{2R},
                 \end{array}
               \right.
\end{equation*} which implies that when $S_1\ne 0$,
\begin{equation*}
    \lim_{R\to +\infty}\frac{1}{R^2}\int_{B_{2R}}S_1dM=+\infty.
\end{equation*}
The proof is complete.
\cqd

\begin{myth}\label{thm3.2}
There is no complete  immersed strongly stable hypersurface
$M^n\rightarrow \mathbb{R}^{n+1}$, $n\geq 3$, with positive constant
scalar curvature and polynomial growth of $1$-volume, that is
$$\lim_{R\rightarrow\infty}\frac{\int _{B_R}S_1dM}
{R^n}<\infty,$$ where $B_R$ is a geodesic ball of radius $R$ of
$M^n$.
\end{myth}
{\bf Proof.} Suppose that $M$ is  a complete  immersed strongly
stable hypersurface $M^n\rightarrow \mathbb{R}^{n+1}$, $n\geq 3$,
with positive constant scalar curvature. From Theorem \ref{thm2.2}, it
suffices to show that the $1$-volume $\int_MS_1dM$ is infinite which
is the part (2) of Lemma \ref{lem3.2}.\\ \cqd

\section {
 Graphs with $S_2=const$ in Euclidean space}
In this section we include some stability properties and estimates
for entire graphs on $\mathbb{R}^n$ which may be known to experts
ant not easy to find a reference. Using these facts we give the
proof of Corollary \ref{coroD}. Let $M^n$ a hypersurface of
$\mathbb{R}^{n+1}$ given by a graph of a function
$u:\mathbb{R}^n\rightarrow \mathbb{R}$ of class
$\mathcal{C}^\infty(\mathbb{R}^n)$. For such hypersurfaces we have:
\begin{myprop}
Let $M^n$ a graph of a function $u:\mathbb{R}^n\rightarrow
\mathbb{R}$ of class $\mathcal{C}^\infty(\mathbb{R}^n)$. Then
\begin{enumerate}
       \item If $S_2=0$  and $S_1$ does not change sign on $M$, then $M^n$ is a stable hypersurface.
    \item If $M$ has $S_2=C> 0$, then $M^n$ is strongly stable.
\end{enumerate}

\end{myprop}
{\bf Proof.} Considerer 
and $f:M\rightarrow
\mathbb{R}$ a $\mathcal{C}^\infty$ function with compact support and let
  $W=\sqrt{1+|\nabla u|^2}$.
In order to calculate $\langle P_1 (\nabla f),\nabla f\rangle$,
write $g=fW$. Thus
\begin{eqnarray*}
  \langle P_1 (\nabla f),\nabla f\rangle &=& \langle P_1 (\nabla (\frac{g}{W})),\nabla (\frac{g}{W})\rangle \\
  &=&\langle P_1 (g\nabla \frac{1}{W}+ \nabla g \frac{1}{W}),g\nabla (\frac{1}{W})+\frac{1}{W}\nabla g\rangle  \\
 &=&\langle g P_1 (\nabla \frac{1}{W})+ \frac{1}{W} P_1(\nabla g ),g\nabla (\frac{1}{W})+\frac{1}{W}\nabla g\rangle  \\
  &=& g^2 \langle P_1 (\nabla \frac{1}{W}), \nabla \frac{1}{W}\rangle+\frac{g}{W} \langle  P_1 (\nabla \frac{1}{W}), \nabla g\rangle\\
  &&+\frac{g}{W} \langle P_1 (\nabla g), \nabla \frac{1}{W}\rangle+\frac{1}{W^2}\langle P_1 (\nabla g), \nabla g\rangle.
\end{eqnarray*}
By using that $P_1$ is selfadjoint, we have:
\begin{equation}\label{graf1}
    \langle P_1 (\nabla f),\nabla f\rangle =g^2 \langle P_1 (\nabla \frac{1}{W}), \nabla \frac{1}{W}\rangle+2\frac{g}{W} \langle  P_1 (\nabla \frac{1}{W}), \nabla g\rangle+\frac{1}{W^2}\langle P_1 (\nabla g), \nabla g\rangle.
\end{equation}

On the other hand, if $\{e_1,...,e_n\}$ is a geodesic frame along
$M$,
\begin{eqnarray*}
  \mathrm{div}(fgP_1(\nabla\frac{1}{W}) &=&\sum_{i=1}^n\langle \nabla_{e_i}(fgP_1(\nabla\frac{1}{W})),e_i\rangle  \\
   &=& \sum_{i=1}^n\langle fg_iP_1(\nabla\frac{1}{W})+f_igP_1(\nabla\frac{1}{W})+fg \nabla_{e_i}(P_1(\nabla\frac{1}{W})),e_i\rangle \\
  &=& \sum_{i=1}^n\{fg_i\langle P_1(\nabla\frac{1}{W}),e_i\rangle+f_ig\langle P_1(\nabla\frac{1}{W}),e_i\rangle+fg \langle\nabla_{e_i}(P_1(\nabla\frac{1}{W})),e_i\rangle\}.
\end{eqnarray*}
 Since $f=\dfrac{g}{W}$, we get
 $$f_i=g_i\frac{1}{W}+g\left(\frac{1}{W}\right)_i,$$
that is,
\begin{eqnarray*}
  gf_i &=& gg_i\frac{1}{W}+g^2\left(\frac{1}{W}\right)_i \\
   &=& fg_i+g^2\left(\frac{1}{W}\right)_i
\end{eqnarray*}
Hence,
\begin{eqnarray*}
   \mathrm{div}(fgP_1(\nabla\frac{1}{W})&=& \sum_{i=1}^n\{fg_i\langle P_1(\nabla\frac{1}{W}),e_i\rangle+(fg_i+g^2\left(\frac{1}{W}\right)_i)\langle P_1(\nabla\frac{1}{W}),e_i\rangle\}+fg L_1(\frac{1}{W})\\
   &=&\sum_{i=1}^n\{2fg_i\langle P_1(\nabla\frac{1}{W}),e_i\rangle+g^2\left(\frac{1}{W}\right)_i\langle P_1(\nabla\frac{1}{W}),e_i\rangle\}+fg L_1(\frac{1}{W})\\
   &=&2f\langle P_1(\nabla\frac{1}{W}),\nabla g\rangle+g^2\langle P_1(\nabla\frac{1}{W}),\nabla (\frac{1}{W})\rangle+fg L_1(\frac{1}{W})\\
   &=&2\frac{g}{W}\langle \nabla\frac{1}{W},P_1(\nabla g)\rangle+g^2\langle P_1(\nabla\frac{1}{W}),\nabla (\frac{1}{W})\rangle+f^2WL_1(\frac{1}{W}).
\end{eqnarray*}
Thus,
\begin{equation}\label{graf2}
    2\frac{g}{W}\langle \nabla\frac{1}{W},P_1(\nabla g)\rangle=
    \mathrm{div}(fgP_1 (\nabla\frac{1}{W}))-g^2\langle P_1(\nabla\frac{1}{W}),\nabla (\frac{1}{W})\rangle-f^2WL_1(\frac{1}{W}).
\end{equation}
Now, by using (\ref{graf2}) into equation (\ref{graf1}), we get
$$
\langle P_1 (\nabla f),\nabla f\rangle =
\mathrm{div}(fgP_1(\nabla\frac{1}{W}))-f^2WL_1(\frac{1}{W})+\frac{1}{W^2}\langle
P_1(\nabla g), \nabla g\rangle.$$

Now, the divergence theorem implies that
\begin{equation*}
    \int_M\langle P_1 (\nabla f),\nabla f\rangle dM= -\int_M f^2WL_1(\frac{1}{W})dM+\int_M\frac{1}{W^2}\langle P_1(\nabla g), \nabla g\rangle dM.
\end{equation*}
Choose the orientation of $M$ in such way that $S_1\geq 0$. Since $S_1^2-|A|^2=2S_2\geq 0$, we obtain that $S_1\geq |A|$. Thus, $\langle P_1(\nabla g), \nabla g\rangle=S_1 |\nabla g|^2 - \langle A\nabla g, \nabla g \rangle \geq (S_1-|A|)|\nabla g|^2\geq 0$, which implies that
\begin{equation}\label{graf3}
    \int_M\langle P_1 (\nabla f),\nabla f\rangle dM\geq -\int_M f^2WL_1(\frac{1}{W})dM.
\end{equation}
When $S_2$ is constant, we will use the following formula proved by Reilly   (see \cite{Re}, Proposition C).
$$L_1(\frac{1}{W})=L_1(\langle N, e_{n+1}\rangle)+(S_1S_{2}-3S_3) \langle N, e_{n+1}\rangle =0,$$
where $N$ is the normal vector of $M$ and $e_{n+1}=(0,...,0, \pm 1),
$ according to our choice of the orientation of $M$.

Thus,
$$  \int_M\langle P_1 (\nabla f),\nabla f\rangle dM\geq -\int_M f^2WL_1(\frac{1}{W})dM=0$$
for all function $f$ with compact support. Hence $M$ is stable if $S_2=0$ and strongly stable in the case $S_2\neq 0$.

\cqd
\begin{myrem}We would like to remark that the operator $L_1$ need not to be elliptic in the above
proof.
\end{myrem}
\begin{myprop}
Let $M^n$ a graph of a function $u:\mathbb{R}^n\rightarrow
\mathbb{R}$ of class $\mathcal{C}^\infty(\mathbb{R}^n)$, with $S_1
\geq 0$. Let $B_R$ be a geodesic ball of radius $R$ in $M$. Then
$$\int_{B_{\theta R}}S_1dM \leq\frac{C(n)}{1-\theta}R^{n},$$
where $C(n)$ and $ \theta$ are constants, with $0<\theta <1$. In
particular, ${\displaystyle \int _MS_1 dM}$ has polynomial growth.
\end{myprop}
{\bf Proof.} Let $f:M\rightarrow \mathbb{R}$ be a be a function in
$\mathcal{C}^\infty_0(M)$, that is a smooth function with compact
support. Observe that $$\mathrm{div}\left (f\frac{\nabla u}{W}\right)=f
\mathrm{div}\left (\frac{\nabla u}{W}\right)+\left\langle \nabla
f,\frac{\nabla u}{W}\right\rangle,$$ where $W=\sqrt{1+|\nabla
u|^2}$. By using the fact that $S_1$ is given by $ \displaystyle{
  S_1 = \mathrm{div}\left( \frac{ \nabla u}{W}\right)}$, we have that

\begin{equation}\label{eqgraf}\int_M f S_1dM=\int_M f\mathrm{div}\left( \frac{\nabla u}{W} \right )
dM= -\int_M \left\langle \nabla f,\frac{\nabla u}{W}\right\rangle
dM.\end{equation}

Now, choose a family of geodesic balls $B_R$ that exhausts $M$. Fix
$\theta$, with $0<\theta\ <1$ and let $f:M\rightarrow \mathbb{R}$ be
a continuous function that is one on $B_{\theta R}$, zero outside
$B_R$ and linear on $B_R\setminus B_{\theta R}$. Therefore, from
equation (\ref{eqgraf}) we obtain
\begin{equation*}
    \int _{B_{\theta R}} S_1 dM \leq \int_{B_R} fS_1 dM \leq \int_{B_R}
    \left \langle \frac{\nabla u}{W}, \nabla f\right \rangle dM.
\end{equation*}
By using Cauchy-Schwarz inequality and the fact that $\dfrac{|\nabla
u|}{W}\leq 1$, it follows that
$$
\int _{B_{\theta R}} S_1 dM\leq \int_{B_{R}}|\nabla f| dM \leq
\int_{B_r\setminus B_{\theta R}} \frac {1}{(1-\theta)R} dM \leq
\frac {1}{(1-\theta)R} \mathrm{vol} (B_R).
$$
We observe that since $M$ is a graph, if $\Omega_R=\{
(x_1,..,x_{n+1})\in \mathbb{R}^{n+1}| -R\leq x_{n+1}\leq R;\;$ $
\sqrt {x_1^2+...+x_n^2}\leq R\} $, then
$$\mathrm{vol} (B_R) \leq \int_{\Omega_R} 1 dx_1...dx_{n+1}=C(n)R^{n+1}.$$
Hence,
$$
\int _{B_{\theta R}} S_1 dM\leq \frac {1}{(1-\theta)R} \mathrm{vol} (B_R)=
\frac{C(n)}{1-\theta} R^n.
$$\\ \cqd

We have the following Corollary of Theorem \ref{thm3.2}
\begin{mycoro}\label{coroD}
Any entire
graph on $\mathbb{R}^n$ with nonnegative constant scalar curvature
must have zero scalar curvature.
\end{mycoro}
{\bf Proof.} Suppose by sake of contradiction that there exist a entire graph with $S_2=const>0$. Such graph is strongly stable and if $S_2> 0$, we get that $S_1^2= |A|^2+ 2S_2>0$, we obtain that $S_1$ does not change sign and we can choose the orientation in such way that $S_1>0$. Thus the graph has polynomial growth of the $1$-volume. Thus we have a contradiction with Theorem \ref{thm3.2}. Thus it follows that $S_2=0$.\\ \cqd

\parbox[t]{3.3in}{
Hil\'ario Alencar \\
Instituto de Matem\'{a}tica\\
Universidade Federal de Alagoas\\
57072-900 Macei\'o-AL, Brazil\\
hilario@mat.ufal.br}
\parbox[t]{3.3in}{
Walcy Santos \\
Instituto de Matem\'{a}tica\\
Universidade Federal do Rio de Janeiro\\
Caixa Postal 68530\\
21941-909, Rio de Janeiro-RJ,
Brazil\\
walcy@im.ufrj.br}
\parbox[t]{3.3in}{
 Detang Zhou\\ Instituto de Matem\' atica\\
Universidade Federal Fluminense\\   24020-140, Niter\'{o}i-RJ Brazil \\
zhou@impa.br}
\end{document}